\documentclass{article}
\usepackage{latexsym}
\usepackage{amsfonts}
\def\Bbb{\mathbb}
\def\bcp{{\Bbb C \Bbb P}}
\def\bea{\begin{eqnarray*}}
\def\eea{\end{eqnarray*}}
\def\mv{\mbox{Vol}}
\newtheorem{main}{Theorem}

\newtheorem{thm}{Theorem}
\newtheorem{prop}{Proposition}
\newtheorem{cor}{Corollary}
\newtheorem{lem}{Lemma}
\newenvironment{proof}{\medskip \noindent
{\bf Proof.}}{\hfill \rule{.5em}{1em}
\\}

\newenvironment{rmk}{\mbox{ }\\{\bf  Remark}\mbox{ }}{
\hfill $\Box$\mbox{}\bigskip}
\begin{document}
\sloppy
\title{Weyl Curvature, Einstein Metrics,\\
and Seiberg-Witten Theory}

\author{Claude LeBrun\thanks{Supported 
in part by  NSF grant DMS-9505744.} 
\\ 
SUNY Stony
 Brook 
  }

\date{March 20, 1998\\
Revised April 14, 1998}
\maketitle

\begin{abstract} We show that solutions of the Seiberg-Witten equations 
lead to non-trivial estimates for the $L^{2}$-norm of the
Weyl curvature of a smooth compact 4-manifold. These 
estimates are then used to derive 
new obstructions to the 
existence of Einstein metrics on smooth compact 4-manifolds
with a non-zero Seiberg-Witten invariant. These 
results considerably refine those previously obtained \cite{lno}
by using scalar-curvature estimates alone.  
 \end{abstract}

\section{Introduction}

A smooth Riemannian metric $g$ is said to be {\em  Einstein} 
if its Ricci curvature $r$ is a constant multiple of 
the metric: 
$$r=\lambda g .$$
Not every smooth compact oriented 4-manifold $M$
admits such a metric. Indeed, 
a well-known necessary condition  is that 
$M$ must satisfy \cite{hit,thorpe,bes} 
the 
Hitchin-Thorpe inequality $2\chi  (M) \geq  3 |\tau (M)|$,
where $\chi$ and $\tau$ respectively denote the signature and
Euler characteristic. This is forced on one by 
the Gauss-Bonnet-like formula
\begin{equation}
(2\chi \pm 3\tau ) (M) 
=\frac{1}{4\pi^{2}}\int_{M}\left(
2|W_{\pm}|^{2}+\frac{s^{2}}{24} -\frac{|\stackrel{\circ}{r}|^{2}}{2}
\right) d\mu ,	
	\label{bonnet}
\end{equation}
where $s$, $\stackrel{\circ}{r}$, $W_{+}$, and $W_{-}$
 respectively denote the scalar, trace-free
Ricci, self-dual Weyl, and anti-self-dual Weyl curvatures of 
an arbitrary Riemannian metric $g$, whose point-wise  
norms $|\cdot |$ and  volume form $d\mu$ also
appear in the integral. The Hitchin-Thorpe inequality 
 follows immediately because Einstein 
 metrics are characterized by the vanishing of 
 $\stackrel{\circ}{r}$, and $\stackrel{\circ}{r}$
 makes the only negative contribution to the integrand.

One could obviously improve this result if  one 
had, say,  non-trivial 
 lower bounds for the the scalar-curvature contribution
 to the integral.
 And indeed, this is 
 precisely what has happened over the course 
 of the past few years.
  For example, simplicial volume \cite{grom,bes} 
  and entropy estimates \cite{bcg2,samba}
lead to new obstructions for certain spaces with 
infinite fundamental group. In another direction, 
the Hitchin-Thorpe argument can  be dramatically
improved upon \cite{lmo,lno} through the use of 
scalar-curvature estimates arising from the Seiberg-Witten
equations  \cite{witten,KM} if, for example, the smooth
4-manifold in question admits a symplectic form \cite{taubes}. 

 The present article will show that
 the Seiberg-Witten equations also give rise to {\em a priori}
 control of the $L^{2}$-norm of $W_{+}$. Our main result in this
 direction is as follows:

\begin{main} \label{rest}
Let $(M,g)$ be a compact oriented Riemannian 4-manifold with
 a non-trivial Seiberg-Witten invariant.  
Let $c_{1}(L)\in H^{2}(M, {\Bbb R})$ be the 
first Chern class of the corresponding spin$^{c}$ structure 
on $M$, and let $c_{1}^{+}\neq 0$ denote its projection into the 
space of $g$-self-dual harmonic 2-forms. Then  
$$\frac{1}{4\pi^{2}} \int_M 
\left( 2|W_{+}|^{2}+\frac{s^2}{24}\right)~d\mu > 
\frac{32}{57}(c_1^+)^2 . $$
 \end{main}

 This leads to new obstructions to the 
 existence of Einstein metrics. Recall that a 
 complex surface $M$ of general type admits a 
 {\em K\"ahler}-Einstein metric only if $M$ 
 cannot be smoothly
 decomposed as a connected sum $X\# k\overline{\bcp}_{2}$, $k > 0$. 
 If one instead wishes to consider Einstein metrics
 which are not necessarily K\"ahler, similar statements can be proved
  if one assumes  that $k$ is sufficiently large. 
 Indeed, Theorem \ref{rest} implies the following  result:

\begin{main} \label{blowup}
Let $X$ be a smooth compact oriented 4-manifold with 
$2\chi + 3\tau > 0$. Assume, moreover, that $X$ has a non-trivial
Seiberg-Witten invariant. Then 
$X\# k \overline{\bcp}_{2}$ does not admit 
Einstein metrics for any $k \geq \frac{25}{57}(2\chi + 3\tau )(X).$ 
\end{main}

This result should be compared to the main result of 
\cite{lno}, where the same conclusion is reached 
for $k \geq \frac{2}{3}(2\chi + 3\tau )(X)$. On the other hand, 
 the Hitchin-Thorpe inequality 
 merely excludes existence when $k\geq (2\chi + 3\tau )(X)$.
 Since $\frac{25}{57} < \frac{4}{9}= (\frac{2}{3})^{2}$, 
 it seems fair to say that Theorem \ref{blowup}
  improves on \cite{lno} by a bigger margin than 
  that by which \cite{lno} 
 improved upon its antecedent.

It should be emphasized that 
the obstructions studied here 
strongly depend on the given smooth structure.
Indeed, Kotschick  \cite{kot} recently pointed out
infinitely many  examples  which 
do not admit Einstein metrics by \cite{lno}, but which are 
nonetheless   homeomorphic to 4-manifolds which {\em do}
admit Einstein metrics. As we shall see in \S \ref{ex}, Theorem \ref{blowup} 
gives rise to a plethora of examples of this
type  which would have been unobtainable by previous methods.

\section{Weyl Estimates}
Let $M$ be a smooth, compact, oriented 4-manifold. Each 
Riemannian metric $g$ on $M$ then determines
a direct sum decomposition
$$H^{2}(M, {\Bbb R}) =
{\cal H}^{+}_{g}\oplus {\cal H}^{-}_{g},$$
where ${\cal H}^{+}_{g}$  (respectively, ${\cal H}^{-}_{g}$)
 consists of those cohomology
classes for which the   harmonic  representative  is 
 self-dual  (respectively,  anti-self-dual).
 The non-negative integer
$b_{+}(M)= \dim {\cal H}^{+}_{g}$  is independent of 
$g$, and we will henceforth always assume it to be positive. 
It is thus  natural to consider the set of metrics $g$
for which ${\cal H}^{+}_{g}=H$ for some 
fixed $b_{+}(M)$-dimensional subspace $H\subset H^{2}(M, {\Bbb R})$;
such metrics will be said to be {\em $H$-adapted}. 
 Assuming there is at least one $H$-adapted metric, we will then 
 say that $H$ is a {\em polarization} of $M$, and \cite{lpm} 
 call the pair 
$(M,H)$  a {\em polarized 4-manifold}. Notice that
the restriction of
the intersection pairing
$$\smile : H^{2}(M, {\Bbb R})\times H^{2}(M, {\Bbb R})\to {\Bbb R}$$
to $H$ 
is then positive definite, and that $H\subset H^{2}$
is   maximal among  subspaces with this property.

Let $c$ be a spin$^{c}$ structure on $M$. Then $c$ determines 
a Hermitian line-bundle $L\to M$ with 
$$c_{1}(L)\equiv w_{2}(M)\bmod 2,$$
and for each metric $g$ we also  have  
rank-2 complex vector bundles
$V_{\pm}\to M$ which formally satisfy 
$$V_{\pm}={\Bbb S}_{\pm}\otimes L^{1/2},$$
where ${\Bbb S}_{\pm}$ are the locally-defined left- and right-handed 
spinor bundles of $g$. Given a polarization $H$ on $M$, we 
will then use $c_{1}^{+}$ to denote the 
orthogonal projection of $c_{1}(L)$ into $H$ with 
respect to the intersection form. If $g$ is a particular
metric with ${\cal H}^{+}_{g}=H$, we will also freely use
$c_{1}(L)$ to denote the $g$-harmonic 2-form representing the 
corresponding 
de Rham class, and use $c_{1}^{+}$ to denote  its self-dual part. 
For example, if a choice of $H$-compatible metric $g$
has already been made,  the number 
$$|c_{1}^{+}|:=\sqrt{(c_{1}^{+})^{2}}$$
may freely be identified with the
L$^{2}$-norm of the self-dual $g$-harmonic form 
denoted by $c_{1}^{+}$.

 For each Riemannian metric $g$,  
the Seiberg-Witten equations \cite{witten} 
\begin{eqnarray} D_{A}\Phi &=&0\label{drc}\\
 F_{A}^+&=&i \sigma(\Phi)\label{sd}\end{eqnarray}
are 
equations 
for an unknown Hermitian connection $A$ on $L$
and an unknown 
smooth section $\Phi$ of $V_+$.
Here the canonical real-quadratic map $\sigma : 
V_{+}\to \Lambda^{+}$ is invariant under parallel transport, and 
satisfies $|\sigma (\Phi )|^{2}= |\Phi |^{4}/8$. 
 If a spin$^c$ structure satisfies 
 $[c_{1}(L)]^{2}=(2\chi + 3\tau )(M)$,  and if 
  $c_{1}^{+}\neq 0$ relative to the polarization $H={\cal H}^{+}_{g}$, 
 then   the  Seiberg-Witten invariant  
$n_c(M , H)$ can be defined \cite{KM,witten}  as the number of solutions,
modulo gauge transformations and  counted with orientations, of a  
generic perturbation 
\begin{eqnarray*} D_{A}\Phi &=&0  \\
 iF^+_A+\sigma (\Phi ) &=& \phi  \end{eqnarray*}
 of  (\ref{drc}--\ref{sd}), 
where $\phi$ is a  smooth self-dual 2-form of small 
$L^{2}$ norm.   More generally, if
 $c$ is  a spin$^c$ structure 
for which  $\ell=[c_1^2(L) - (2\chi + 3\tau)(M)]/4$ is non-negative 
and even, one can 
define  \cite{t2}  the perturbed Seiberg-Witten invariant 
 $n_c(M , H)$  to be $\int_Z \eta^{\ell/2}$, where
$Z$ is the $\ell$-dimensional smooth compact moduli space of solutions of 
a generic small perturbation of the equations,
and $\eta  \in H^2(Z)$ is the first Chern class of
the based moduli space, considered as an $S^1$-bundle over $Z$.
For other values of $c_1^2(L)$, one simply sets $n_{c}(M,H)=0$ as a
matter of definition. 
 For our purposes, the  point 
is simply that when $c_{1}^{+}$ and  $n_{c}(M,H)$
are both non-zero,
(\ref{drc}--\ref{sd}) must 
have a solution with $\Phi\not\equiv 0$
 for every $H$-adapted metric $g$. Moreover, these equations 
imply the Weitzenb\"ock formula
\begin{equation}\label{wb}
 4\nabla_A^*\nabla_A \Phi +  s \Phi +|\Phi|^2\Phi  =0.
\end{equation}

\begin{thm} \label{est}
Let $(M,H)$ be a  polarized smooth compact  
oriented  4-manifold, and let $c$ be a spin$^{c}$ structure
on $M$  for which the  
 Seiberg-Witten invariant  $n_{c}(M,H)$ 
is non-zero; let $c_1(L)\in H^2(M, {\Bbb R})$ denote 
the anti-canonical class of $c$, and 
  let $c_1^+\neq 0$ be its orthogonal projection to
$H$ with respect to the intersection form. 
 Then 
every 
$H$-adapted Riemannian metric $g$ satisfies
$$\left(\int |W_{+}|^{2} ~d\mu\right)^{{1/2}}
\geq (4-3{\beta}^{1/2} )\frac{2\pi}{\sqrt{3}} |c_{1}^{+}| ,$$
where 
\begin{equation}
\beta =\frac{ \int s^{2}d\mu}{32\pi^{2}(c_{1}^{+})^{2}} \geq 1.	
	\label{alpha}
\end{equation}
Moreover, equality can  occur only if $\beta =1$. 
\end{thm}

\begin{proof}
Recall that every every conformal class on $M$ contains \cite{aubin,lp,sch}
a Yamabe 
metric, and that such a metric   minimizes $\int s^{2}d\mu$ 
in its conformal class \cite{bcg1,lno}. 
Now all such  metrics have constant scalar 
curvature. Since $\int |W_{+}|^{2}d\mu$ is conformally invariant, 
the form of the desired inequality therefore allows us to 
 assume henceforth that
$g$ has  constant  scalar curvature. But because 
$(c^{+}_{1})^{2}\neq 0$ and 
the Seiberg-Witten invariant $n_c(M, H)$
 is assumed to be non-zero, $g$ cannot \cite{lpm} 
have $s\geq 0$. We may thus assume henceforth that $g$
has constant scalar curvature $s < 0$. 

 Now because $n_c(M, H)\neq 0$ by assumption, 
 there must exist an
irreducible solution 
of  (\ref{drc}) and  (\ref{sd}). 
But the Weitzenb\"ock formula (\ref{wb}) tells
us that 
\begin{eqnarray*}
8\int_M  |\nabla_A\Phi|^2 ~d\mu	 & = &2 \int  (-s-|\Phi|^2)|\Phi|^2~d\mu   \\
	 & = &
	    \int  (-s-|\Phi|^2)(-s+|\Phi|^2)~d\mu -  
	    \int  (-s-|\Phi|^2)^2~d\mu\\
	 & = & \|s\|_{2}^{2} - \|\Phi \|_{4}^{4}
	 -  \|(|s|-|\Phi |^{2}) \|^{2}_{2} \\
	 & \leq  & \|s\|_{2}^{2} -  \|\Phi \|^{4}_{4}
	 -(\|s\|_{2}-\|\Phi \|_{4}^{2})^{2} \\
	 & = & 2(\|s\|_{2}-\|\Phi \|_{4}^{2})\|\Phi \|_{4}^{2}
\end{eqnarray*}
and hence
$$\frac{\|\nabla_A\Phi\|_{2}^{2}}{\|\Phi \|_{4}^{2}} 
\leq \frac{1}{4}(\|s\|_{2}-4\sqrt{2}\pi |c_{1}^{+}|)=
(\beta^{1/2}-1)\sqrt{2}\pi |c_{1}^{+}| .
$$
 Here we have used the observation that
$|\Phi |^{4}= 8|F^+_A|^{2}$, as is implied by (\ref{sd}),
together with the fact that $2\pi c_{1}^{+}$ is the harmonic part of
$F^+_A$.  

 Since we have assumed that $s$ is a negative constant, 
the Weitzenb\"ock formula (\ref{wb}) also implies \cite{KM} 
the $C^{0}$ estimate 
$$| \Phi |^{2} \leq |s|.$$
Since (\ref{drc}) implies that 
$$|\nabla F^+_A|^{2}\leq  
\frac{1}{2}|\Phi |^{2} 
|\nabla_A\Phi |^{2}$$
it now follows that 
$$\frac{\|\nabla F_{A}^{+}\|_{2}^{2}}{\| F_{A}^{+} \|_{2}} 
\leq \sqrt{2} |s| \frac{\|\nabla_A\Phi\|_{2}^{2}}{\|\Phi \|_{4}^{2}} 
\leq  |s|
(\beta^{1/2}-1)2\pi |c_{1}^{+}| .$$

Now any self-dual  2-form $\varphi$ on any oriented 4-manifold satisfies 
\cite{bourg} the 
Weitzenb\"ock formula 
$$(d+d^{*})^{2}\varphi = \nabla^{*}\nabla \varphi - 2W^{+}(\varphi , 
\cdot ) + \frac{s}{3} \varphi,$$
where $W^{+}$ is the self-dual Weyl tensor. It follows that 
$$\int_{M}(-W^{+})(\varphi , \varphi ) \geq 
\int_{M}(-\frac{s}{6})|\varphi |^{2}~d\mu -
 \frac{1}{2}\int_{M} |\nabla \varphi |^{2}
 ~d\mu .$$ On the other hand, 
 $W^{+}$ is a trace-free quadratic form on $\Lambda^{+}$, so
 that $|W^{+}| |\varphi |^{2} \geq \sqrt{\frac{3}{2}}
 (-W^{+})(\varphi , \varphi )$. Again assuming that the scalar 
 curvature $s$ of $g$ is a negative constant,
 we therefore have
 $$
 \left(\int |W^{+}|^{2}d\mu \right)^{1/2}
  \left(\int |\varphi|^{4}d\mu \right)^{1/2}
  \geq \frac{|s|}{2\sqrt{2}} \left[ \frac{1}{\sqrt{3}}
  \|\varphi\|_{2}^{2}-\frac{\sqrt{3}}{|s|}
  \|\nabla\varphi \|_{2}^{2}
  \right]
 $$
 and hence, assuming that  $\varphi \not\equiv 0$,
 we have 
$$
\left(\int |W^{+}|^{2}d\mu \right)^{1/2}
  \geq \sqrt{{\frac{s^{2}\|\varphi\|_{2}^{2}}{8\|\varphi\|_{4}^{4}}}}
  \left[ \frac{1}{\sqrt{3}}
  \|\varphi\|_{2}-\frac{\sqrt{3}}{|s|}
 \frac{ \|\nabla\varphi \|_{2}^{2}}{\|\varphi\|_{2}}
  \right].$$

We now 
 apply this  to  $\varphi =  F^+_A$.
Because $|F^+_A|^{2}= \frac{1}{8}|\Phi|^{4}\leq \frac{s^2}{8}$,
we have 
$$
{\frac{s^{2}\|F^+_A\|_{2}^{2}}{8\|F^+_A\|_{4}^{4}}}
=\frac{(s^{2}/8)\int |F^+_A|^{2}d\mu}{\int |F^+_A|^{4}d\mu}\geq 1 .
$$
It follows that 
  \begin{eqnarray*} 
 \left(\int |W^{+}|^{2}d\mu \right)^{1/2}	 & \geq  & 
 \frac{1}{\sqrt{3}}
  \|F^+_A\|_{2}-\frac{\sqrt{3}}{|s|}
 \frac{ \|\nabla F^+_A  \|_{2}^{2}}{\|F^+_A\|_{2}}\\
 	 & \geq &  \frac{1}{\sqrt{3}}\|F^+_A\|_{2}-\frac{\sqrt{3}}{|s|}  |s|
(\beta^{1/2}-1) 2\pi |c_{1}^{+}|   \\
  	 & \geq  & \frac{ 2\pi }{\sqrt{3}}  |c_{1}^{+}| 
- \sqrt{3}   (\beta^{1/2}-1) 2\pi |c_{1}^{+}|  \\
  	 & =  & \left[ 1- 3 (\beta^{1/2}-1) \right] 
  	 \frac{2\pi}{\sqrt{3}} |c_{1}^{+}|   \\
  	 & = & \left( 4- 3 \beta^{1/2} \right)   
  	 \frac{2\pi}{\sqrt{3}} |c_{1}^{+}| . 
  \end{eqnarray*}
 If equality holds, moreover, 
 $F_{A}^{+}$ must be harmonic and $| F_{A}^{+} |^{2}\equiv s^{2}/8$,
 so that  $\beta =1$, as claimed.   
\end{proof}

This immediately implies a new characterization of 
constant-scalar-curvature 
K\"ahler metrics; cf. \cite{lpm}.  

\begin{cor}\label{sharp}
Let $(M,H)$ be a  polarized smooth compact  
oriented  4-manifold, and let $c$ be a spin$^{c}$ structure
on $M$. Assume, moreover, that  $c_{1}^{+}\neq 0$ and that the   
 Seiberg-Witten invariant  $n_{c}(M,H)$ is non-zero. 
 Then 
every 
$H$-adapted Riemannian metric $g$ on $M$ satisfies
$$\frac{1}{4\pi^{2}} \int_M 
\left( \frac{1}{3}|W_{+}|^{2}+\frac{s^2}{24}\right)~d\mu \geq 
\frac{4}{9}(c_1^+)^2 , $$
with equality iff $g$  is   K\"ahler  (for a
$c$-compatible complex structure)  
and has constant negative scalar curvature. 
 \end{cor}
\begin{proof}
Let us set 
$$\alpha = \frac{3}{2}(\beta -1) = \frac{3}{2}\left(\frac{\int 
s^{2}d\mu}{32\pi^{2}(c_{1}^{+})^{2}} -1 \right) \geq 0, $$
so that 
$$\beta^{{1/2}}\leq 1 + \frac{1}{3}\alpha . $$
Theorem \ref{est} then tells us that 
$$\left( \int_M  |W_{+}|^{2} ~d\mu \right)^{{1/2}}
\geq (1-\alpha ) \frac{2\pi}{\sqrt{3}} |c_{1}^{+}|. $$
But $x \geq 1-h \Longrightarrow x^{2}\geq 1-2h ,$ 
since the parabola $y=x^{2}$ sits above its tangent line
at $x=1$.  
We therefore have  
$$
\int_M  |W_{+}|^{2} ~d\mu 
\geq (1-2\alpha ) \frac{4\pi^{2}}{3}
(c_{1}^{+})^{2}=
4\pi^{2}\frac{4}{3}
(c_{1}^{+})^{2}- \frac{1}{8}\int s^{2}d\mu , 
$$
which is to say that 
$$\frac{1}{4\pi^{2}}\int_{M}\left(
\frac{1}{3}|W_{+}|^{2}
+\frac{s^{2}}{24}\right) ~d\mu
\geq \frac{4}{9}(c_{1}^{+})^{2}.$$

If equality is achieved, the metric is a Yamabe minimizer, and 
 so has constant scalar curvature; moreover, 
$\nabla F^+_A=0$, and since $F^+_A\not\equiv 0$, it follows that 
 the metric is K\"ahler. 
 Conversely, the Seiberg-Witten invariant is  non-zero for
 a K\"ahler surface with $c_{1}\cdot [\omega ] < 0$,
 where $[\omega ]$ is the
  K\"ahler class; 
  and 
 since any K\"ahler metric satisfies $|W_{+}|^{2}=s^{2}/24$,
  $\int s~d\mu = 4\pi c_{1}\cdot [\omega ]$, and 
  $\int d\mu =[\omega ]^{2}/2$, the inequality is precisely saturated by
  K\"ahler metrics of constant negative scalar curvature.
\end{proof}

\begin{rmk}
If $c_{1}^{+}=0$ and $b_{+} =1$, the Seiberg-Witten invariant
is ill-defined. However, 
 the relevant inequality 
$$\frac{1}{4\pi^{2}} \int_M 
\left( \frac{1}{3}|W_{+}|^{2}+\frac{s^2}{24}\right)~d\mu \geq 
\frac{4}{9}(c_1^+)^2  =0$$
has become a triviality in this case. 
Moreover,  equality  occurs
 precisely when
the metric is scalar-flat and anti-self-dual, and, since
 we have 
 assumed that $b_{+} \neq 0$,  the Weitzenb\"ock formula
 for 2-forms then shows \cite{laf,lsd} that 
 any such metric is scalar-flat K\"ahler. 
 Much the same conclusion therefore holds in this 
 case, albeit for slightly different  reasons.
\end{rmk}

In another direction, recall that Taubes \cite{tasd} has shown 
that for any smooth compact orientable $X^{4}$, there is an integer
$k_{0}$ such that $M=X\# k\overline{\bcp}_{2}$ admits metrics
with $W_{+}=0$ provided that $k \geq k_{0}$. In particular, there
are many anti-self-dual 4-manifolds with non-trivial Seiberg-Witten 
invariants. For such manifolds, we get an interesting 
scalar-curvature estimate.

\begin{cor}\label{asdm}
Let $(M,g)$ be a compact {\em anti-self-dual} 
4-manifold with a non-zero Seiberg-Witten invariant.  Then 
 $$\frac{1}{32\pi^{2}}\int_{M}s^{2}d\mu > 
 \frac{16}{9}(c_{1}^{+})^{2},$$
 where $c_{1}^{+}$ is again the self-dual part of the 
 first Chern class of the relevant spin$^{c}$ structure.
 \end{cor}
 
\begin{proof}
By Theorem \ref{est}, one must have $\beta > (\frac{4}{3})^{2}$ if 
$W_{+}\equiv 0$.
\end{proof}

One may use this as a vanishing theorem. For example, it 
immediately 
implies that the Seiberg-Witten invariant must vanish for all those
spin$^{c}$ structure on a hyperbolic 4-manifold which satisfy
 $(c_{1}^{+})^{2} \geq 
\frac{32}{3}\chi$. One might  guess that all 
the  Seiberg-Witten invariants of 
a hyperbolic 4-manifold must vanish, but there is very
 little hard evidence 
either  for or against such a conjecture.

We now come to our main technical result, which, while certainly
not sharp in the above sense, will yield better results 
in many interesting contexts:

 \setcounter{main}{0}
\begin{main}\label{keen}
Let $(M,H,c)$ be as above. 
 Then 
every 
$H$-adapted Riemannian metric $g$ satisfies
$$\frac{1}{4\pi^{2}} \int_M 
\left( 2|W_{+}|^{2}+\frac{s^2}{24}\right)~d\mu > 
\frac{32}{57}(c_1^+)^2 . $$
 \end{main}

\begin{proof}
Our definition (\ref{alpha}) of $\beta$ has been chosen so that
$$
\frac{1}{4\pi^{2}}\int_{M}\frac{s^{2}}{24}d\mu= \beta
\frac{(c_{1}^{+})^{2}}{3}.$$
If $\beta > \frac{16}{9}$, we therefore have 
$$\frac{1}{4\pi^{2}} \int_M 
\left( 2|W_{+}|^{2}+\frac{s^2}{24}\right)~d\mu > \frac{16}{27}
(c_{1}^{+})^{2} > \frac{32}{57}(c_1^+)^2  $$
for trivial reasons. We may therefore assume henceforth that $\beta 
\in [1,\frac{16}{9}]$.

This assumption, however, guarantees that both sides of the inequality
$$\left(\int |W_{+}|^{2} ~d\mu\right)^{{1/2}}
\geq (4-3\beta^{{1/2}} )\frac{2\pi}{\sqrt{3}} |c_{1}^{+}| $$
are non-negative. It then follows that 
$$
\frac{1}{4\pi^{2}}
\int |W_{+}|^{2} ~d\mu \geq (4-3\beta^{{1/2}})^{2} \frac{(c_{1}^{+})^{2}}{3}.
$$
Hence
\begin{eqnarray*}
\frac{1}{4\pi^{2}} \int_M 
\left( 2|W_{+}|^{2}+\frac{s^2}{24}\right)~d\mu 	 & \geq  &
 [\beta + 2(4-3\beta^{1/2} 
)^{2}]\frac{(c_{1}^{+})^{2}}{3}  \\
	 & = & \left[19 (\beta^{1/2}-\frac{24}{19})^{2}+ \frac{32}{19}\right]  
	  \frac{(c_{1}^{+})^{2}}{3}
\end{eqnarray*}
with equality only if $\beta=1$. Hence 
$$
\frac{1}{4\pi^{2}} \int_M 
\left( 2|W_{+}|^{2}+\frac{s^2}{24}\right)~d\mu 
>   \frac{32}{57} (c_{1}^{+})^{2}
$$
whenever the Seiberg-Witten invariant is non-zero.
\end{proof}

\begin{rmk} Kronheimer \cite{K} recently showed that certain
 4-manifolds with vanishing Seiberg-Witten invariants nonetheless
 have the remarkable property that 
 there is a
  solution of the Seiberg-Witten equations
  (with  fixed spin$^{c}$ structure)
 for each and every metric. The present results manifestly also apply
 to such manifolds, since the above proofs stem directly from 
 structure 
 of the equations rather than from formal properties of 
 the invariant. For our present purposes, however, this phenomenon 
 does not yet seem to have any interesting ramifications;  in 
 particular, 
all of  Kronheimer's examples have $2\chi + 3\tau \leq 0$, and 
indeed 
can actually  be collapsed to 
 zero volume while keeping $s$ and $|W_{+}|$  bounded. 
\end{rmk}

\section{Einstein Metrics}

Up until now, we have been discussing the Seiberg-Witten invariants
of a polarized 4-manifold $(M,H)$. These are always  well-defined for
any spin$^{c}$ structure with $c_{1}^{+}\neq 0$. When $b_{+}(M)> 1$,
moreover, they are even independent of the choice of $H$, since any two 
generic 
choices of polarization may be joined by a path for which $c_{1}^{+}$
is always non-zero. When $b_{+}(M)=1$, however, a different feature 
emerges: the intersection form is of Lorentz type, and 
the value of the invariant depends on whether $c_{1}^{+}$ is 
`future pointing' or `past pointing' with respect to a given 
time-orientation for  $H^{2}(M)$. On the other hand, if our 
manifold has $2\chi + 3\tau > 0$, $c_{1}(L)$ is automatically
time-like, and  only one of these possibilities actually occurs. 
Thus it makes perfectly good sense to speak of the the Seiberg-Witten
invariant of a 4-manifold with $b_{+}=1$ as long as 
$2\chi + 3\tau > 0$. Indeed, the same reasoning even applies if
$2\chi + 3\tau =0$, provided that $c_{1}(L)$ is not a torsion class. 

The following observation \cite{lno,kot} is the work-horse 
to which our estimates will be harnessed:

\begin{lem} \label{who} 
Let $X$ be a smooth compact oriented 4-manifold with
$2\chi + 3\tau > 0$. Assume, moreover, that some 
Seiberg-Witten invariant of $X$ is non-trivial. 
Let $k$ be any natural number, and let  $H$ be any  polarization of 
$M= X\# k \overline{\bcp}_{2}$. Then there is a spin$^{c}$
structure $c$ on $M$ such that $n_{c}(M,H)\neq 0$ and 
$$(c_{1}^{+})^{2}\geq (2\chi + 3\tau ) (X).$$  
\end{lem}
\begin{proof}
Let $c_{1}(X)$ denote the first Chern class of a
spin$^{c}$ structure on $X$ for which the Seiberg-Witten invariant
is non-zero, and notice that $(c_{1}(X))^{2}\geq
 (2\chi + 3\tau ) (X) > 0$, because the relevant Seiberg-Witten 
 moduli space must have non-negative virtual dimension.  
 Pull $c_{1}(X)$ back to $M= X\# k \overline{\bcp}_{2}$ via
 the canonical collapsing map, and, by a standard  abuse of 
 notation, let us also use $c_{1}(X)$ to denote the 
 pulled-back class. Thus, with respect to our given polarization,
 $$([c_{1}(X)]^{+})^{2}\geq (c_{1}(X))^{2}\geq
 (2\chi + 3\tau ) (X) > 0.$$
 Now choose  generators
 $E_{1}, \ldots , E_{k}$ for the  pull-backs to $M$ of
  the 
 $k$ relevant copies of $H^{2}(\overline{\bcp}_{2}, {\Bbb Z})$  
 so that 
 $$[c_{1}(X)]^{+}\cdot E_{j}\geq 0, ~~~j = 1, \ldots , k.$$
 Then \cite{fs} there is a spin$^{c}$ structure on $M$ with
 $n_{c}(M,H)\neq 0$ and 
 $$c_{1}(L)= c_{1}(X) + \sum_{j=1}^{k}  E_{j}.$$
 But one then has
 \begin{eqnarray*}
 	(c_{1}^{+})^{2} & = & \left([c_{1}(X)]^{+} + 
 	\sum_{j=1}^{k} 	E_{j}^{+}\right)^{2}  \\
 	 & = & ([c_{1}(X)]^{+})^{2}+ 2 \sum_{j=1}^{k} [c_{1}(X)]^{+}\cdot  E_{j}
 + (\sum_{j=1}^{k}  E_{j}^{+})^{2}  \\
 	 & \geq  & ([c_{1}(X)]^{+})^{2}  \\
 	 & \geq  & (2\chi + 3\tau ) (X), 
 \end{eqnarray*}
 exactly as claimed. 
\end{proof}

Our main result  now follows.

\setcounter{main}{1}
\begin{main} 
 Let $X$ be a smooth compact oriented 4-manifold with 
$2\chi + 3\tau > 0$. Assume, moreover, that $X$ has a non-trivial
Seiberg-Witten invariant. Then 
$X\# k \overline{\bcp}_{2}$ does not admit an
Einstein metric if $k \geq \frac{25}{57}(2\chi + 3\tau )(X).$ 
\end{main}
\begin{proof}
  For any  Einstein 
metric $g$ on $M$, equation (\ref{bonnet}) and Theorem \ref{keen}
tell us that 
$$ (2\chi + 3\tau )(M)=\frac{1}{4\pi^2}\int_M \left(2|W_+|^2 +  
\frac{s^2}{24}\right)d\mu > \frac{32}{57} (c_{1}^{+})^{2} $$
for any spin$^{c}$ structure with $n_{c}(M,H)\neq 0$, where
$H={\cal H}^{+}_{g}$. 
But Lemma \ref{who} now asserts that $M=X\# k \overline{\bcp}_{2}$
has such a spin$^{c}$ structure with $(c_{1}^{+})^{2} \geq (2\chi + 
3\tau)(X)$. Thus
$$
(2\chi + 3\tau)(X) - k = (2\chi + 3\tau)(M) > \frac{32}{57}
(2\chi + 3\tau)(X) ,$$
and hence 
$$ k < \frac{25}{57}
(2\chi + 3\tau)(X),$$
assuming that $M$ admits an Einstein metric. 
The result therefore follows by contraposition. 
\end{proof}

\begin{cor}   \label{nein}
Let $X$ be a minimal complex  surface
of general type, 
and let $M= X\#k\overline{\bcp}_2$ be obtained from
$X$ by blowing up $k$ points.
If $k \geq  \frac{25}{57}  c_1^2(X)$,
then $M$ does not admit  Einstein metrics. 
 \end{cor}

One might instead ask whether $M=X\# k \overline{\bcp}_{2}$
admits {\em anti-self-dual} Einstein metrics, since Taubes'
theorem  \cite{tasd}
tells us that anti-self-dual (but non-Einstein)
metrics exist when  $k$ is very large.  Using Corollary
\ref{asdm} instead of Theorem \ref{keen}, Lemma \ref{who} 
 implies a slightly better estimate by essentially the same argument:

\begin{prop}
 Let $X$ be a smooth compact oriented 4-manifold with 
$2\chi + 3\tau > 0$ and a non-trivial
Seiberg-Witten invariant. Then 
$X\# k \overline{\bcp}_{2}$ cannot admit anti-self-dual 
Einstein metrics if $k \geq  \frac{11}{27}(2\chi + 3\tau )(X).$ 
\end{prop}

\section{Examples}\label{ex}

We will now examine some  specific   new examples of 4-manifolds
without Einstein metrics given to us by Corollary \ref{nein}.

Let us begin by considering
 a non-singular complex hypersurface $X_{\ell}$ of 
degree $\ell > 4$ in $\bcp_{3}$. This minimal complex surface of 
general type has $c_{1}^{2}=\ell (\ell -4 )^{2}$ and 
$p_{g}= {\ell -1 \choose 3}.$ If we blow up $X_{\ell}$ at
at $k \geq \frac{25}{57} \ell (\ell -4 )^{2}$ points,
the result is a complex surface 
$M = X_{\ell} \# k \overline{\bcp}_{2}$ which is not 
diffeomorphic to any Einstein manifold. In particular,
$X_{9}\# 117 \overline{\bcp}_{2}$ does not admit
Einstein metrics. But this  
 complex surface 
has 
$c_{1}^{2}= 108$ and $p_{g}= 56$, exactly like the 
the  double-branched cover of
$\bcp_{1}\times \bcp_{1}$ ramified over a 
curve of bidegree $(6,58)$. This so-called Horikawa surface \cite{hori},
like $X_{9}\# 117 \overline{\bcp}_{2}$, is simply
connected, and both surfaces have $\tau = 
c_{1}^{2}-8(1+p_{g})= -348\not\equiv 0 \bmod 16$,
so Freedman's classification theorem
\cite{freedman} tells us that both are homeomorphic
to 
$$113 \bcp_{2}\# 461 \overline\bcp_{2}. 
$$
Our Horikawa surface, however, has
ample canonical line bundle, and so admits 
a K\"ahler-Einstein metric by the Aubin/Yau theorem \cite{aubin0,yau}. 
Thus, while
$X_{9}\# 117 \overline{\bcp}_{2}$  does not
admit  Einstein metrics, it is nonetheless 
homeomorphic to an Einstein manifold. Similar examples
may be obtained by blowing up the hypersurface
$X_{\ell}$ for any $\ell \geq 9$. The reader  might
wish to compare this with the results of 
\cite{kot}, where the above class of examples was 
observed to be beyond the capabilities of 
the weaker non-existence result  of \cite{lno}. 

 	 One may   construct more delicate examples
 	by blowing up  branched double covers	
 	 $Y_{m}$ of $\bcp_{2}$,	with ramification locus  
 	 a	smooth curve of	degree $2m$, $m> 4$. These are surfaces
 	 of general type with $c_{1}^{2}=2(m-3)^{2}$, and 
 	 Corollary	\ref{nein} yields new examples
 	 of 4-manifolds without Einstein metrics by considering
 	 their	blow-ups. In	particular,	
 	 $Y_{27}\#	506\overline{\bcp}_{2}$, does not 
 	 admit	Einstein metrics, since  
 	 $506 >	\frac{25}{57} 2\cdot 24^{2}$. However, this
 	 simply connected  complex surface
 	  has $c_{1}^{2}= 646$	and	
 	 $p_{g}= 325$, and it follows that it is 
 	 homeomorphic to a Horikawa surface \cite{hori} --- e.g. 
 	 the double branched cover of 
 	 $\bcp_{2}\# \overline\bcp_{2}$,  ramified over
 	 the proper transform of a curve of degree 
 	 $330$ in $\bcp_{2}$ which is non-singular 
 	 except for $324$ self-crossings at the blown up point. 
 	 In the same way, one can show that for all $m \geq 27$,
 	  the double planes  $Y_{m}$ have blow-ups which
 	  are homeomorphic to Einstein manifolds, but do 
 	  not themselves admit Einstein metrics.

\section{Minimality and Minimal Volumes}

Let us now turn to  a  discussion of
 minimal volume problems \cite{grom,bcg2,lno}.  Given 
 a compact smooth 4-manifold $M$,   let 
 ${\cal M}_{|r|}(M)$ and ${\cal M}_s(M)$
 respectively denote
the set of metrics on $M$ for which $3g\geq r\geq -3g$ and 
$s\geq -12$. Then we may define  minimal volume invariants 
$$\mv_{|r|}(M):= \inf_{g\in {\cal M}_{|r|}}\int_Md\mu_g $$ 
 $$\mv_s(M):= \inf_{g\in {\cal M}_s}\int_Md\mu_g  ,$$
and our conventions have been
chosen so that $\mv_{|r|}(M)\geq \mv_s(M)$ tautologically. 
If $M$ is a complex surface of general type and 
$X$ is its minimal model, it was observed in \cite{lno} that 
$$\mv_s(M) = \mv_s(X) = \mv_{|r|}(X) = 
\frac{2}{9}  \pi^{2}c_{1}^{2}(X).$$ 
If $M$ is non-minimal, however, we will now see that
 $\mv_{|r|}(M) > \mv_{s}(M)$. 

\begin{lem}
Let $(M,H)$ be a  polarized smooth compact  
oriented  4-manifold, and let $c$ be a spin$^{c}$ structure
on $M$. Assume, moreover, that  $c_{1}^{+}\neq 0$ and that the   
 Seiberg-Witten invariant  $n_{c}(M,H)$ is non-zero. 
 Then 
every 
$H$-adapted Riemannian metric $g$ on $M$ satisfies
$$\frac{1}{8\pi^{2}} \int_M |r_{g}|^{2}d\mu \geq 
\frac{8}{5}(c_1^+)^2-\frac{3}{5}(2\chi + 3\tau)(M) , $$
with equality iff the metric is K\"ahler-Einstein. 
\label{which}
 \end{lem}
 \begin{proof}
 The Gauss-Bonnet formula (\ref{bonnet}) tells us that 
 	 \begin{eqnarray*}
	 	\frac{1}{8\pi^{2}} \int_M |r|^{2}d\mu & = & 
 \frac{1}{8\pi^{2}}  \int_M \left(\frac{s^{2}}{4}+ 
 |\stackrel{\circ}{r}|^{2} \right)  d\mu  \\
	 	 & = & \frac{1}{\pi^{2}}\int_{M}
	 \left(\frac{s^{2}}{24}+\frac{1}{2}|W_{+}|^{2} \right)
	 d\mu - (2\chi + 3\tau)(M) .
	 \end{eqnarray*}
	 On the other hand, the same formula tells us that  
	 $$
 \frac{1}{4\pi^{2}}\int_{M}
	 \left(\frac{s^{2}}{24}+2|W_{+}|^{2} \right)
	 d\mu \geq (2\chi + 3\tau)(M) ,$$
	 with equality iff the metric is Einstein. 
	 Meanwhile, 
  Corollary
 \ref{sharp} asserts that 
  $$\frac{1}{4\pi^{2}}\int_{M}
	 \left(\frac{s_{g}^{2}}{24}+\frac{1}{3}|W_{+}|_{g}^{2} \right)
	 d\mu_{g} \geq  \frac{4}{9}(c_1^+)^2 , $$
	 with equality iff the metric is constant-scalar-curvature K\"ahler.
	 Taking a convex combination of these inequalities,  
	 we therefore have 
	 \begin{eqnarray*}
	 \frac{1}{4\pi^{2}}\int_{M}
	 \left(\frac{s_{g}^{2}}{24}+\frac{1}{2}|W_{+}|_{g}^{2} \right)
	 d\mu_{g} 	 & \geq  & \frac{9}{10}\cdot \frac{4}{9}(c_1^+)^2 
	 + \frac{1}{10}(2\chi + 3\tau)(M)  \\
	 	 & = & \frac{2}{5}(c_1^+)^2 
	 + \frac{1}{10}(2\chi + 3\tau)(M). 
	 \end{eqnarray*}
	Hence 
\begin{eqnarray*}
	 \frac{1}{8\pi^{2}} \int_M |r_{g}|^{2}d\mu_{g}& \geq & 
	 \frac{8}{5}(c_1^+)^2 
	 + \frac{2}{5}(2\chi + 3\tau)(M)-(2\chi + 3\tau)(M)  \\
	 & = & 
\frac{8}{5}(c_1^+)^2-\frac{3}{5}(2\chi + 3\tau)(M) ,
\end{eqnarray*}
and, as claimed,  equality  holds precisely in the K\"ahler-Einstein 
case. 
 \end{proof}
 
 Now if $X$ is a minimal complex surface of general type, 
 and if $M=X\# k\overline{\bcp}_{2}$ is one of its blow-ups, 
 Lemmata \ref{who} and \ref{which} assert that every metric 
 on $M$ satisfies 
 \begin{eqnarray*}
 	\frac{1}{8\pi^{2}} \int_M |r_{g}|^{2}d\mu_{g} &  \geq  & 
 	 \frac{8}{5}c_{1}^{2}(X) -\frac{3}{5}\left[c_{1}^{2}(X)-k
 \right]   \\
 	 & = & c_{1}^{2}(X) + \frac{3}{5}k . 
 \end{eqnarray*}
 Since any metric with $-3g\leq r\leq 3g$ satisfies $|r|^{2}\leq 36$,
 this shows that every  $g\in {\cal M}_{|r|}(M)$ has total volume
 \begin{eqnarray*}
 	 \int_{M} d\mu_{g} & \geq  & 
 	 \frac{8\pi^{2}}{36}[c_{1}^{2}(X) + \frac{3}{5}k]  \\
 	 & = & \frac{2}{9}\pi^{2} c_{1}^{2}(X)  + \frac{2}{15}\pi^{2}k,  
 \end{eqnarray*}
 so that 
 \begin{eqnarray*}
 \mv_{|r|}(M)	 & \geq & 
  \frac{2}{9}\pi^{2} c_{1}^{2}(X)  + \frac{2}{15}\pi^{2}k  \\
 	 & = & \mv_{s}(M)+ \frac{2}{15}\pi^{2}k .
 \end{eqnarray*}
 This proves

\begin{thm} \label{vol}
Let $X$ be a    complex surface of general type. Then 
$$\mv_{|r|}(M) = \mv_{s}(M) \Longleftrightarrow \mbox{ $M$ is minimal.}$$
\end{thm}

The reader should note that, while Lemma \ref{which}
provides an effective means of proving Theorem
\ref{vol}, the relevant estimate is actually quite weak in practice.
For instance, an argument parallel to the proof of Theorem \ref{keen}
shows that 	 $$\frac{1}{4\pi^{2}}\int_{M}
	 \left(\frac{s_{g}^{2}}{24}+\frac{1}{2}|W_{+}|_{g}^{2} \right)
	 d\mu_{g} > \frac{16}{33} (c_{1}^{+})^{2}$$
	 whenever the Seiberg-Witten invariant is non-zero.
 On a non-minimal surface $M=X\# k\overline{\bcp}_{2}$ 
 of general type, every metric 
	 therefore satisfies
$$
\frac{1}{8\pi^{2}}	\int_{M}|r|^{2}_{g}d\mu_{g}  >   
\frac{31}{33}c_{1}^{2}(X) +k ,$$
so that 
 $$\mv_{|r|}(M)\geq \frac{31}{33}\mv_{s}(M) + \frac{2}{9}\pi^{2}k ,
 $$
 and this, of course, is a stronger estimate for all but the 
 smallest   values of $k$.
 On the other hand, there is every reason to expect that
 the   actual value of
 $\mv_{|r|}(M)$ might in fact  be considerably larger than 
 indicated by either of these estimates; after all, 
 the definition of   $\mv_{|r|}(M)$   essentially 
 involves the 
  $C^{0}$ norm, rather than its $L^{2}$ norm, of $r$.

 In light of the above discussion, it is natural to
 introduce some new smooth invariants of a smooth oriented 4-manifold 
 $M$. 
For any real parameter $\varepsilon \geq 0$, let us set 
  $${\cal I}_{\varepsilon}(M) := \inf_{g} \frac{1}{4\pi^{2}}\int_{M} 
  \left(\frac{s_{g}^{2}}{24}
  + \varepsilon |W_{+}|_{g}^{2}\right) d\mu_{g}.$$
  Thus, for example, the above discussion of $\mv_{{|r|}}$
  hinged on estimating ${\cal I}_{\frac{1}{2}}(M)$. 
  On the other hand, 
  $96\pi^{2}{\cal I}_{0}(M)= [\min(  Y(M), 0)]^{2}$, where $Y(M)$ is the
  Yamabe invariant (sigma constant) of $M$. When $\varepsilon$ is small,
  Corollary \ref{sharp} actually allows us to compute this invariant
  exactly for most complex surfaces:
  
  \begin{thm}
  Let $(M,J)$ be a complex surface of Kodaira dimension $\geq 0$,
  and let $(X,\tilde{J})$ be its minimal model. Then 
  $${\cal I}_{\varepsilon}(M)= \frac{1+\varepsilon }{3}(2\chi + 3\tau)(X)
  $$ 
  for any $\varepsilon\in [0,\frac{1}{3}]$. 
  \end{thm}
  \begin{proof}
  If $M$ has Kodaira dimension $0$ or $1$, it 
  collapses \cite{lky} 
  with $s$ and $|W_{+}|$ bounded, and hence 
  ${\cal I}_{\varepsilon}(M)=0$ for all 
  $\varepsilon\geq 0$.
  Since $(2\chi + 3\tau)(X)=c_{1}^{2}(X)=0$ in this case, the result
  is established in Kodaira dimensions $0$ and $1$. 
  
  If, on the other hand, $M$ is a surface of general type, we can 
  proceed as in \cite{lno}, constructing metrics on $M$ 
  by gluing   gravitational instantons and Burns metrics
  onto  a K\"ahler-Einstein orbifold, 
  while  changing the $L^{2}$-norms of $s$ and $W_{+}$ as little
  as we like. In this way, we see that 
  $${\cal I}_{\varepsilon}(M)\leq \frac{1+\varepsilon}{3}c_{1}^{2}(X) =  
  \frac{1+\varepsilon}{3}(2\chi +  3\tau)(X) ~~\forall \varepsilon \geq 0.$$
 Yet Corollary \ref{sharp} and Lemma \ref{who} 
 tell us that any metric on $M$ satisfies
 $$\frac{1}{4\pi^{2}}\int_{M} 
  \left(\frac{s_{g}^{2}}{24}
  + \frac{1}{3} |W_{+}|_{g}^{2}\right) d\mu_{g} \geq 
  \frac{4}{9} c_{1}^{2}(X), $$
  establishing the opposite inequalities when $\varepsilon= 1/3$.
  However, the inequality 
 $$\frac{1}{4\pi^{2}}\int_{M} 
 \frac{s_{g}^{2}}{24}
  d\mu_{g} \geq 
  \frac{1}{3} c_{1}^{2}(X) $$
  was 
  already proved in \cite{lno}. Taking 
  convex combinations, we therefore have 
 $${\cal I}_{\varepsilon}(M)\geq \frac{1+\varepsilon}{3}c_{1}^{2}(X)$$ 
  for all $\varepsilon\in [0, \frac{1}{3}]$,  and the result follows.   
  \end{proof}
  
When the Kodaira dimension is zero, one can still read off
${\cal I}_{\varepsilon}$  for most complex surfaces. 
Indeed, if $M$ is the underlying
4-manifold of a ruled surface of genus $\geq 2$, then
\cite{lrs,kp} one can find
scalar-flat anti-self-dual metrics on $M$, so that 
${\cal I}_{\varepsilon}(M)=0$.
If the base has genus 1, $M$ is diffeomorphic
 to an elliptic surface, so that it collapses with 
 bounded $s$ and $W_{+}$ as in \cite{lky},
 and hence has ${\cal I}_{\varepsilon}=0$. The same 
 argument applies to the 
 rational elliptic surface 
 $\bcp_{2}\# 9 \overline{\bcp}_{2}$ and all of its blow-ups. 
 
 Now the above  discussion covers all  
   complex algebraic surfaces except for  
 the rational surfaces with 
 $2\chi + 3\tau > 0$. On the other hand, equation (\ref{bonnet})
  implies that 
 $${\cal I}_{\varepsilon}(M) \geq \min ( 1, \frac{\varepsilon}{2}) 
 (2\chi + 3\tau ) (M)$$
 We
 therefore have the following result:
 
 \begin{prop}
 Let $(M^{4},J)$ be a complex algebraic surface, and let 
 $\varepsilon$ be any positive constant. 
 Then ${\cal I}_{\varepsilon}(M)> 0$ iff 
 \begin{itemize}
 \item $(M,J)$ is of general type; or
 \item $(M,J)$ is a  deformation of a del Pezzo surface. 
 \end{itemize}
 \end{prop}

\vfill
\noindent
{\bf Acknowledgment.}
This work grew out of calculations
carried out while the author was a visiting scholar
at Harvard University. He would therefore like to
thank the members of the Harvard Mathematics Department, and especially
Cliff Taubes, for their hospitality during his stay. 

\pagebreak

\end{document}